\documentclass{IEEEtran}
\IEEEoverridecommandlockouts

\usepackage{cite}
\usepackage{amsmath,amssymb,amsfonts}
\usepackage[capitalise]{cleveref}
\usepackage[]{footmisc}
\usepackage{xcolor}
\usepackage{algorithm}
\usepackage{algorithmic}
\usepackage[algo2e]{algorithm2e} 
\usepackage{graphicx}
\usepackage{textcomp}
\usepackage{psfrag}
\usepackage{epsfig}
\usepackage{subcaption}
\usepackage{float}
\usepackage{xcolor}

\usepackage{booktabs}
\usepackage{multicol}
\usepackage{arydshln}
\usepackage{multirow}
\usepackage{wrapfig}
\usepackage[numbers,sort&compress]{natbib}

\newtheorem{mypro}{Proposition}
\newtheorem{myass}{Assumption}
\newtheorem{myrem}{Remark}
\newtheorem{mylem}{Lemma}
\newtheorem{mythr}{Theorem}

\newtheorem{myprob}{Problem}

\SetKwRepeat{Repeat}{Repeat}{Until}
\SetKw{Return}{Return}

\begin{document}

\title{\LARGE Data-Driven Structured Control for Continuous-Time LTI Systems}

\author{Zhaohua~Yang,~Yuxing~Zhong,~and~Ling~Shi,~\IEEEmembership{Fellow,~IEEE}
	\thanks{\rm Z. Yang, Y. Zhong, and L. Shi are with the Department of Electronic and Computer Engineering, Hong Kong University of Science and
		Technology, Clear Water Bay, Hong Kong SAR (email: zyangcr@connect.ust.hk; yuxing.zhong@connect.ust.hk; eesling@ust.hk). Ling Shi is also with the Department of Chemical and Biological Engineering, Hong Kong University of Science and Technology, Clear Water Bay, Hong Kong SAR.}
	\thanks{}
  }
\maketitle

\begin{abstract}
    This paper addresses the data-driven structured controller design problem for continuous-time linear time-invariant (LTI) systems. We consider three control objectives, including stabilization, $H_2$ performance, and $H_\infty$ performance. Using the collected data, we construct a minimal matrix ellipsoid that contains all
    admissible system matrices. We propose some linearization techniques that enable us to incorporate the structural constraint on the controller, which motivates an iterative algorithm for each control objective. Finally, we provide some numerical examples to demonstrate the effectiveness of the proposed methods.
\end{abstract}


\IEEEpeerreviewmaketitle
\section{Introduction}
Networked control systems (NCSs) have found various applications in modern industry. In such systems, controllers establish a communication network between sensors and actuators to achieve desired control objectives. The structure of the communication network is thus reflected in the sparsity pattern of the controller. In many practical applications, the communication network is subject to structural constraints due to the inherent physical and structural properties of the system. This motivates the study of structured controller design, which aims to design a controller with a specific sparsity pattern that satisfies the desired control objectives \cite{jovanovic2016controller}. The structured controller design problem is known to be NP-hard in general \cite{blondel1997np}. However, various methods have been proposed to tackle this problem. For instance, Lin et al. \cite{lin2011augmented} proposed an augmented Lagrangian method to design the optimal structured $H_2$ controller. Fardad and Jovanovi{\'c} \cite{fardad2014design} developed an iterative linear matrix inequality (ILMI) algorithm to design the optimal structured $H_2$ controller. More recently, Ferrante et al. \cite{ferrante2019design} proposed a sufficient condition that allows the design of structured stabilizing controllers within a convex subset of the original nonconvex feasibility region. They then extended this method to tackle the structured $H_\infty$ controller design problems \cite{ferrante2020lmi}. These works assume the system model is known \textit{a priori}. However, in many practical applications, the system model is unknown or difficult to obtain. This motivates our study of data-driven controller design.

Data-driven control has received significant attention in recent years due to its ability to design controllers directly from data without requiring an explicit system model. In a seminal work, De Persis and Tesi \cite{de2019formulas} developed a data-driven framework to design stabilizing controllers, considering both noise-free and noisy data. Under the common assumption that the data are bounded, van Waarde et al. \cite{van2020noisy} proposed the matrix S-lemma, a powerful tool that establishes the connection among quadratic matrix inequalities (QMIs). They then applied this lemma to the robust stabilization problem and the $H_2$ ($H_\infty$) control problem. Additionally, Bisoffi et al. \cite{bisoffi2022data} introduced the concept of matrix ellipsoid and solved the robust stabilization problem based on Petersen's lemma. For data-driven structured control, Miller et al. \cite{miller2025data} extended the work in \cite{ferrante2019design} to the data-driven setting and designed structured $H_2$ controllers. Yang et al. \cite{yang2025data} improved the results in \cite{miller2025data} by developing iterative algorithms, and they further extended the results to the $H_\infty$ control case. However, these works focus on discrete-time systems, and the continuous-time case remains unexplored. Moreover, the analyses in \cite{miller2025data,yang2025data} are based on the matrix S-lemma, which is not directly applicable to the continuous-time case. Therefore, new methods are needed to address the data-driven structured control problem for continuous-time systems.

Our main contributions are summarized as follows:   
\begin{enumerate}
    \item Unlike existing literature, where the system model is perfectly known, we, for the first time, investigate structured control for continuous-time unknown LTI systems using data. We consider three control objectives, including stabilization, $H_2$ performance, and $H_\infty$ performance.
    \item Using the collected data, we first construct a minimal matrix ellipsoid containing all possible system matrices that can generate the data. Then, we introduce novel linearization techniques to relax the original nonconvex problems into convex SDPs, which motivates an effective algorithm for each control objective.
    \item We compare our methods with the one proposed in \cite{jovanovic2016controller} and demonstrate the advantages of our approach.
\end{enumerate}

This paper is organized as follows. In \cref{Preliminaries}, we introduce the system model, present some preliminary results, and formulate the problem of interest. In \cref{Main Results}, we characterize the uncertainty set of possible system matrices using data and present our main results on data-driven structured controller synthesis for stabilization, $H_2$ control, and $H_\infty$ control. In \cref{Simulations}, we provide numerical experiments to illustrate the effectiveness of the proposed methods. Finally, we conclude the paper in \cref{Conclusion}.

\textit{Notations}: Let $\mathbb{R}^{n\times m}$ and $\mathbb{R}^n$ denote the sets of real matrices of size $n\times m$ and real column vectors of length $n$. 
For a matrix $X$, let $X^\top$ and $\mathrm{Tr}(X)$ denote its transpose and trace. Let $X_{ij}$ ($[X]_{ij}$) denote the entry of $X$ at the $i$-th row and $j$-th column. The expression $X\succ0$ ($X\succeq 0$) means that $X$ is positive definite (semidefinite). Let $\circ$ denote the elementwise product. Let $I$, $0$, and $\mathbf{1}$ denote the identity matrix, zero matrix, and all-one matrix, with dimension labeled if necessary. Let $\mathcal{H}\left(X\right) \triangleq X + X^\top$, $\mathcal{G}(X) \triangleq X^\top X$ for brevity. We sometimes abbreviate $\left[\begin{smallmatrix}A & B \\ B^\top & C\end{smallmatrix}\right]$ as  $\left[\begin{smallmatrix}A & * \\ B^\top & C\end{smallmatrix}\right]$ to save space.

\section{Preliminaries}\label{Preliminaries}
\subsection{System Modeling}
In this paper, we consider the following continuous-time LTI system:
\begin{equation}\label{basic system}
    \begin{aligned}
    \dot{x}(t) &= Ax(t) + Bu(t) + Gd(t),\\
    y(t) &= Cx(t) + Du(t) + Hd(t),
    \end{aligned}
\end{equation}
where $x(t)\in \mathbb{R}^{n_x}$ is the state, $u(t)\in \mathbb{R}^{n_u}$ is the control input, $d(t)\in \mathbb{R}^{n_d}$ is the disturbance input, and $y(t)\in \mathbb{R}^{n_y}$ is the controlled performance output. We make the standard assumption that the pair $(A,B)$ is stabilizable. The system \eqref{basic system} is regulated by a static state-feedback controller of the form:
\begin{equation}\label{state feedback}
    u(t) = Kx(t),
\end{equation}
where $K\in \mathbb{R}^{n_u\times n_x}$ is the feedback gain matrix to be designed. Thus, the closed-loop system is given by:
\begin{equation}\label{closed-loop system}
    \begin{aligned}
    \dot{x}(t) &= (A+BK)x(t) + Gd(t),\\
    y(t) &= (C+DK)x(t) + Hd(t).
    \end{aligned}
\end{equation}
The closed-loop matrices $A+BK$ and $C+DK$ may be denoted by $A_K$ and $C_K$ in the rest of the paper. The transfer function from the disturbance input $d(t)$ to the controlled output $y(t)$ is given by:
\begin{equation}
    T_{yd}(s) = (C+DK)(sI - (A+BK))^{-1}G + H.
\end{equation}
We denote the $H_2$ norm and $H_\infty$ norm of the transfer function $T_{yd}(s)$ as $\|T_{yd}(s)\|_2$ and $\|T_{yd}(s)\|_\infty$, respectively. To facilitate the presentation, we provide the following well-known results regarding $H_2$ and $H_\infty$ control design. Throughout this paper, we call a matrix (Hurwitz) stable if all its eigenvalues have strictly negative real parts.
\begin{mylem}[\cite{scherer2002multiobjective}]\label{stabiling lemma}
    The closed-loop matrix $A+BK$ is stable if and only if there exists a positive definite matrix $P\succ0$ such that
    \begin{equation}\label{stabiling eq}
        (A+BK)^\top P + P(A+BK) \prec 0.
    \end{equation}
\end{mylem}
\begin{mylem}[\cite{scherer2002multiobjective}]\label{H2 lemma}
    The closed-loop matrix $A+BK$ is stable and the $H_2$ norm $\|T_{yd}(s)\|_2 < \gamma$ if and only if there exists a positive definite matrix $P\succ0$ such that
    \begin{subequations}\label{H2 eq}
        \begin{align}
        &\begin{bmatrix}
            (A+BK)^\top P + P(A+BK) &(C+DK)^\top\\
            C+DK & -I
        \end{bmatrix}\prec 0,\label{H2 eq1}\\
        &G^\top PG \preceq Z, \mathrm{Tr}(Z) < \gamma^2.\label{H2 eq2}
        \end{align}
    \end{subequations}
\end{mylem}

\begin{mylem}[\cite{scherer2002multiobjective}]\label{Hinf lemma}
    The closed-loop matrix $A+BK$ is stable and the $H_\infty$ norm $\|T_{yd}(s)\|_\infty < \gamma$ if and only if there exists a positive definite matrix $P\succ0$ such that
    \begin{equation}\label{Hinf eq}
        \begin{bmatrix}
            \mathcal{H}\left(P(A+BK)\right) &PG &(C+DK)^\top\\
            G^\top P & -\gamma I & H^\top\\
            C+DK & H & -\gamma I
        \end{bmatrix}\prec 0.
    \end{equation}
\end{mylem}
\subsection{Structured Controller}
The controller establishes a communication network between the input layer and the state layer, where each nonzero entry $K_{ij} \neq 0$ indicates that the $i$-th input can access the $j$-th state variable. The optimal controller typically yields a dense $K$, requiring each local controller to access the full state information. However, this control strategy may not be feasible in practice due to inherent physical, communication, or structural constraints in the system. This motivates the design of a structured controller with specific sparsity patterns. The structural constraint can be expressed as $K\in S$, where $S$ is a linear subspace of $\mathbb{R}^{n_u\times n_x}$. The corresponding structural identity matrix is defined as:
\begin{equation}
    [I_{S}]_{ij} \triangleq \begin{cases}
    0, &\text{if}\ K_{ij}=0\ \text{is required},\\
    1, &\text{if}\ K_{ij}\ \text{is a free variable}.
    \end{cases}
\end{equation}
We then define its complement as $I_{S^c} = \mathbf{1}_{n_u\times n_x} - I_{S}$. Consequently, the structural constraint $K\in S$ can be equivalently 
expressed as $K\circ I_{S^c} = 0$.

\subsection{Problem Formulation}
In this paper, we consider the scenario where the ground truth system matrices $(A,B)$ are unknown. Instead, we have access only to measured data consisting of input, state, and state derivative samples collected from the system \eqref{basic system} in the presence of an unknown disturbance $d(t)$. Specifically, we collect $T$ samples at discrete time instants $t_0, t_1, \ldots, t_{T-1}$, where $t_0 = 0$ and $t_{i+1} - t_i = T_s$ for all $i = 0, 1, \ldots, T-2$. The collected data are organized into the following matrices:
\begin{equation}
\begin{aligned}
    X_1  &\triangleq \begin{bmatrix}\dot{x}(t_0) &\dot{x}(t_1) &\dots &\dot{x}(t_{T-1})\end{bmatrix} \in \mathbb{R}^{n_x\times T},\\
    X_0 &\triangleq \begin{bmatrix}x(t_0) &x(t_1) &\dots &x(t_{T-1})\end{bmatrix} \in \mathbb{R}^{n_x\times T},\\
    U_0 &\triangleq \begin{bmatrix}u(t_0) &u(t_1) &\dots &u(t_{T-1})\end{bmatrix} \in \mathbb{R}^{n_u\times T}.
\end{aligned}
\end{equation}
In addition, we make the following assumption regarding the disturbance signal $d(t)$.
\begin{myass}\label{noise assumption}
    The disturbance signal $d(t)$ is bounded, i.e., $\|d(t)\|_2 \leq \epsilon$ for all $t>0$ and some $\epsilon>0$.
\end{myass}

Since the ground truth system matrices $(A,B)$ are unknown, in the rest of this paper, unless explicitly specified as the ground truth, we use $(A,B)$ to denote all possible system matrices that can generate the collected data with $d(t)$ satisfying \cref{noise assumption}. These possible system matrices $(A,B)$ can be expressed as a set, which we denote as $\Sigma$. We will derive its mathematical expression in \cref{Mathematical Expression of set Phi}. 

In this paper, we aim to design a state-feedback controller $u=Kx$ subject to an \textit{a priori} structural constraint $K \in S$. Specifically, we consider stabilization and $H_2$ ($H_\infty$) performance. Since we lack knowledge of the ground-truth system matrices $(A,B)$ and only have access to the uncertainty set $\Sigma$, the designed controller should be robust for all $(A,B)\in \Sigma$. In \cref{Mathematical Expression of set Phi}, we clarify that $\Sigma$ is difficult to characterize, and we instead derive a suitable superset $\Phi \supseteq \Sigma$ for subsequent analysis. The problems are formally formulated as follows:
\begin{myprob}\label{problem 1}
    Using the collected data $(X_0,U_0,X_1)$, design a structured controller $K\in S$ such that the closed-loop matrix $A+B K$ is stable for all $(A,B)\in \Phi$.
\end{myprob}
\begin{myprob}\label{problem 2}
    Using the collected data $(X_0,U_0,X_1)$, design a structured controller $K\in S$ to minimize $\gamma$ such that $\|T_{yd}(s)\|_2 < \gamma$ for all $(A,B)\in \Phi$.
\end{myprob}
\begin{myprob}\label{problem 3}
    Using the collected data $(X_0,U_0,X_1)$, design a structured controller $K\in S$ to minimize $\gamma$ such that $\|T_{yd}(s)\|_\infty < \gamma$ for all $(A,B)\in \Phi$.
\end{myprob}
\begin{myrem}
    The problem formulation is conservative as $\Sigma \subseteq  \Phi$. Therefore, the optimal $\gamma^*$ returned by \cref{problem 2} (\cref{problem 3}) is an upper bound on the maximal $H_2$ ($H_\infty$) norm over all $(A,B)\in \Sigma$. We call this upper bound the $H_2$ ($H_\infty$) performance bound.
\end{myrem}

\section{Main Results}\label{Main Results}
\subsection{Mathematical Expression of $\Sigma$ and $\Phi$}\label{Mathematical Expression of set Phi}
In this section, we derive the mathematical expression of the uncertainty set $\Sigma$ and find a suitable superset $\Phi \supseteq \Sigma$. 

\cref{noise assumption} implies that $Gd(t_i) (Gd(t_i))^\top \preceq \epsilon^2 GG^\top$, which can be reformulated as
\begin{equation}
    \begin{bmatrix}
        I\\
        (Gd(t_i))^\top
    \end{bmatrix}^\top 
    \begin{bmatrix}
        -\epsilon^2 GG^\top &0\\
        0 &1
    \end{bmatrix}
    \begin{bmatrix}
        I\\
        (Gd(t_i))^\top
    \end{bmatrix}\preceq 0.
\end{equation}
Combining this with the system dynamics
\begin{equation}
    \dot{x}(t_i) = Ax(t_i) + Bu(t_i) + Gd(t_i),
\end{equation}
we have
\begin{equation}
    \begin{bmatrix}
        I\\
        A^\top\\
        B^\top
    \end{bmatrix}^\top
    \begin{bmatrix}
        \setlength{\dashlinegap}{0.8pt}
        \begin{array}{c:c}
        c_i &b_i^\top\\
        \hdashline
        b_i &a_i
    \end{array}
        \end{bmatrix}
        \begin{bmatrix}
        I\\
        A^\top\\
        B^\top
    \end{bmatrix}\preceq 0,
\end{equation}
where
\begin{equation*}
    \scalebox{0.97}{$
    \begin{bmatrix}
        \setlength{\dashlinegap}{0.8pt}
        \begin{array}{c:c}
        c_i &b_i^\top\\
        \hdashline
        b_i &a_i
        \end{array}
        \end{bmatrix} \triangleq   \begin{bmatrix}
        \setlength{\dashlinegap}{0.8pt}
        \begin{array}{c:c}
        \dot{x}(t_i)\dot{x}(t_i)^\top-\epsilon^2 GG^\top &-\dot{x}(t_i)\begin{bmatrix}
            x(t_i)\\
            u(t_i)
        \end{bmatrix}^\top\\
        \\[-10pt]
        \hdashline
        \\[-10pt]
        -\begin{bmatrix}
            x(t_i)\\
            u(t_i)
        \end{bmatrix}\dot{x}(t_i)^\top &\begin{bmatrix}
            x(t_i)\\
            u(t_i)
        \end{bmatrix}\begin{bmatrix}
            x(t_i)\\
            u(t_i)
        \end{bmatrix}^\top
        \end{array}
        \end{bmatrix}.
        $}
\end{equation*}
Thus, we can define the set of all possible system matrices $(A,B)$ that can generate the $(x(t_i),u(t_i),\dot{x}(t_i))$ as
\begin{equation}
    \Sigma_i \triangleq \left\{ (A,B) \Bigg|     \begin{bmatrix}
        I\\
        A^\top\\
        B^\top
    \end{bmatrix}^\top
    \begin{bmatrix}
        \setlength{\dashlinegap}{0.8pt}
        \begin{array}{c:c}
        c_i &b_i^\top\\
        \hdashline
        b_i &a_i
        \end{array}
        \end{bmatrix}
    \begin{bmatrix}
        I\\
        A^\top\\
        B^\top
    \end{bmatrix}\preceq 0\right\}.
\end{equation}
Consequently, the set $\Sigma$ of $(A,B)$ pairs that are consistent with the full data set $(X_0,U_0,X_1)$ is given by the intersection of all individual sets  $\Sigma_i, i=0,1,\ldots,T-1$:
\begin{equation}
    \Sigma = \bigcap_{i=0}^{T-1} \Sigma_i.
\end{equation}
It is worth noting that each set $\Sigma_i$ is unbounded, and their intersection $\Sigma$ is challenging to express in a compact form. To address this issue, we adopt the matrix ellipsoid framework from \cite{bisoffi2022data}, which provides an analytical method to encapsulate all possible system matrices $(A,B)$ within a matrix ellipsoid. The matrix ellipsoid for $(A,B)$ is typically expressed as:

\begin{equation}\label{first matrix ellipsoid expression}
    \left\{ (A,B)=Z^\top \Bigg|    (Z-\boldsymbol{\delta})^\top\mathbf{A}(Z-\boldsymbol{\delta})\preceq \mathbf{Q}\right\},
\end{equation}
where $\mathbf{A}\succ 0, \mathbf{Q}\succeq 0$. This can be equivalently expressed as
\begin{equation}\label{second matrix ellipsoid expression}
    \left\{ (A,B) \Bigg|     \begin{bmatrix}
        I\\
        A^\top\\
        B^\top
    \end{bmatrix}^\top
    \begin{bmatrix}
        \setlength{\dashlinegap}{0.8pt}
        \begin{array}{c:c}
        \mathbf{C} &\mathbf{B}^\top\\
        \hdashline
        \\[-10pt]
        \mathbf{B} &\mathbf{A}
        \end{array}
        \end{bmatrix}
    \begin{bmatrix}
        I\\
        A^\top\\
        B^\top
    \end{bmatrix}\preceq 0\right\},
\end{equation}
where $\mathbf{A}\succ 0, \mathbf{B}^\top \mathbf{A}^{-1}\mathbf{B}-\mathbf{C}\succeq 0$. From \eqref{second matrix ellipsoid expression}, we can easily compute $\boldsymbol{\delta}$ and $\mathbf{Q}$ in \eqref{first matrix ellipsoid expression} as $\boldsymbol{\delta} = -\mathbf{A}^{-1}\mathbf{B}, \mathbf{Q} = \mathbf{B}^\top \mathbf{A}^{-1}\mathbf{B}-\mathbf{C}$. We aim to find the matrix ellipsoid of smallest volume containing $\Sigma$, which is achieved by minimizing $-\log\det(\mathbf{A})$. We fix $\mathbf{Q}=I$, i.e., $\mathbf{C} = \mathbf{B}^\top\mathbf{A}^{-1}\mathbf{B}-I$, to avoid variable redundancy. Following \cite{bisoffi2022data}, the minimal ellipsoid can be obtained by solving the following optimization problem:
\begin{subequations}\label{minimal ellipsoid}
\begin{align}
    &\min_{\mathbf{A},\mathbf{B},\{\theta\}_{i=0,1,\ldots,T-1}} \quad  -\log\det(\mathbf{A}) \qquad \mathrm{s.t.}\\
        &  \begin{bmatrix}
            \setlength{\dashlinegap}{0.8pt}
            \begin{array}{cc:c}
            -I &\mathbf{B}^\top &\mathbf{B}^\top\\
            \mathbf{B} &\mathbf{A} &0\\
            \hdashline
            \mathbf{B} &0 &-\mathbf{A}
            \end{array}
            \end{bmatrix} - 
            \begin{bmatrix}
            \setlength{\dashlinegap}{0.8pt}
            \begin{array}{c:c}
            \sum_{i=0}^{T-1} \theta_i \begin{bmatrix}c_i &b_i^\top\\ b_i &a_i\end{bmatrix} &0 \\
            \\[-10pt]
            \hdashline
            0 &0 
            \end{array}
            \end{bmatrix} \preceq 0\\
        & \mathbf{A}\succ 0,\theta_i \geq 0 , i = 1,\dots, T-1.
\end{align}
\end{subequations}
The problem is an SDP and thus can be easily solved. Let $\mathbf{A}^*,\mathbf{B}^*$ be the solution to Problem \eqref{minimal ellipsoid}. The minimal ellipsoid can be expressed by
\begin{equation}
    \Phi = \left\{ (A,B)=Z^\top \Bigg|    (Z-\boldsymbol{\delta}^*)^\top\mathbf{A}^*(Z-\boldsymbol{\delta}^*)\preceq I\right\},
\end{equation}
where $\boldsymbol{\delta}^* = -{\mathbf{A}^*}^{-1}\mathbf{B}^*$. We now introduce a critical alternative expression of $\Phi$, which will be used in the subsequent controller design.
\begin{mypro}[\cite{bisoffi2022data}]\label{critical alternative}
\begin{equation}
    \Phi = \left\{  (\boldsymbol{\delta}^*+{\mathbf{A}^*}^{-\frac{1}{2}}\Gamma)^\top \;\middle|\; \|\Gamma\|\leq 1 \right\}.
\end{equation}
\end{mypro}
We are now ready to solve \cref{problem 1}, \cref{problem 2}, and \cref{problem 3}, which will be presented in the rest of \cref{Main Results}.
\subsection{Stabilizing Control}\label{Stabilizing Control}
Before considering the structural constraint, it is natural to check whether an unstructured stabilizing controller exists. Replacing $\begin{bmatrix}
    A &B\end{bmatrix}$ in \eqref{stabiling eq} using \cref{critical alternative}, we obtain
\begin{equation}
    \mathcal{H}\left(P{\boldsymbol{\delta}^*}^\top \begin{bmatrix}       I\\
        K
    \end{bmatrix}\right)
    +
    \mathcal{H}\left(P\Gamma^\top ({\mathbf{A}^*}^{-\frac{1}{2}})^\top \begin{bmatrix}
        I\\
        K
    \end{bmatrix}\right)\prec 0,
\end{equation}
which holds for all $\|\Gamma\|\leq 1$. This statement is equivalent to
\begin{equation}\label{stabilizing temp1}
    \mathcal{H}\left(P{\boldsymbol{\delta}^*}^\top \begin{bmatrix}       I\\
        K
    \end{bmatrix}\right) + PP + \begin{bmatrix}
        I\\
        K
    \end{bmatrix}^\top {\mathbf{A}^*}^{-1} \begin{bmatrix}
        I\\
        K
    \end{bmatrix} \prec 0
\end{equation}
by applying the strict Petersen's lemma \cite{bisoffi2022data}. After changing variables $X = P^{-1}$ and multiplying both sides by $X$, we obtain 
\begin{equation}
    \mathcal{H}\left({\boldsymbol{\delta}^*}^\top \begin{bmatrix}       I\\
        K
    \end{bmatrix}X\right) + I + X\begin{bmatrix}
        I\\
        K
    \end{bmatrix}^\top {\mathbf{A}^*}^{-1} \begin{bmatrix}
        I\\
        K
    \end{bmatrix}X \prec 0.
\end{equation}
Let $Y=KX$. Using the Schur complement, we can obtain an unstructured stabilizing controller by solving the following feasibility problem:
\begin{subequations}\label{data driven unstructured stabilizing}
    \begin{align}
        &\mathrm{find} \quad  X,Y \qquad \mathrm{s.t.}\\
        &\begin{bmatrix}
            {\boldsymbol{\delta}^*}^\top\begin{bmatrix}
                X\\
                Y
            \end{bmatrix}+\begin{bmatrix}
                X\\
                Y
            \end{bmatrix}^\top\boldsymbol{\delta}^*  + I & \begin{bmatrix}
                X\\
                Y
            \end{bmatrix}^\top \\
            \\[-10pt]
            \begin{bmatrix}
                X\\
                Y
            \end{bmatrix} & -\mathbf{A}^*
        \end{bmatrix} \prec 0\\
        &X \succ 0.
    \end{align}
\end{subequations}
Then, a stabilizing controller gain is computed by $\Tilde{K} = \Tilde{Y}\Tilde{X}^{-1}$, where $(\Tilde{X},\Tilde{Y})$ is a feasible solution to \eqref{data driven unstructured stabilizing}. However, if we impose the structural constraint $K\in S$, the change of variable $Y=KX$ is not applicable, as it leads to an intractable constraint $Y X^{-1} \in S$. To address this issue, we directly work with \eqref{stabilizing temp1}. The Schur complement leads to
\begin{equation}\label{stabilizing temp2}
    \begin{bmatrix}
        P{\boldsymbol{\delta}^*}^\top\begin{bmatrix}
        I\\
        K
    \end{bmatrix}+\begin{bmatrix}
        I\\
        K
    \end{bmatrix}^\top\boldsymbol{\delta}^*P  & \begin{bmatrix}
        I\\
        K
    \end{bmatrix}^\top & P\\
    \\[-10pt]
    \begin{bmatrix}
        I\\
        K
    \end{bmatrix} & -\mathbf{A}^* & 0\\
    P &0 &-I
    \end{bmatrix} \prec 0
\end{equation}

The main challenge now lies in handling the bilinear term involving $P$ and $K$ in the top-left block. We first decompose it as
\begin{equation}
\begin{aligned}
    &P{\boldsymbol{\delta}^*}^\top\begin{bmatrix}
    I\\
    K
\end{bmatrix}+\begin{bmatrix}
    I\\
    K
\end{bmatrix}^\top\boldsymbol{\delta}^*P = \mathcal{H}\left(P{\boldsymbol{\delta}^*}^\top\begin{bmatrix}
    I\\
    K
\end{bmatrix}\right) \\
= &\frac{1}{2}\mathcal{G}\left(\boldsymbol{\delta}^*P+\begin{bmatrix}
    I\\
    K
\end{bmatrix}\right)
-\frac{1}{2}\mathcal{G}\left(\boldsymbol{\delta}^*P-\begin{bmatrix}
    I\\
    K
\end{bmatrix}\right)
\end{aligned}
\end{equation}
Next, we focus on the second term. For given matrices $\Tilde{K}$, $\Tilde{P}$,
\begin{equation}
\begin{aligned}
&\mathcal{G}\left(\boldsymbol{\delta}^*P-\begin{bmatrix}
    I\\
    K
\end{bmatrix}\right)\\
=&\mathcal{G}\left(\boldsymbol{\delta}^*\Tilde{P}-\begin{bmatrix}
    I\\
    \Tilde{K}
\end{bmatrix} + \boldsymbol{\delta}^*(P-\Tilde{P})- \begin{bmatrix}0 \\K-\Tilde{K}\end{bmatrix}\right)\\
=&\mathcal{G}\left(\boldsymbol{\delta}^*\Tilde{P}-\begin{bmatrix}
    I\\
    \Tilde{K}
\end{bmatrix}\right)+ \mathcal{G}\left(\boldsymbol{\delta}^*(P-\Tilde{P})- \begin{bmatrix}0 \\K-\Tilde{K}\end{bmatrix}\right) \\
+& \mathcal{H}\left(\left(\boldsymbol{\delta}^*\Tilde{P}-\begin{bmatrix}
    I\\
    \Tilde{K}
\end{bmatrix}\right)^\top\left(\boldsymbol{\delta}^*(P-\Tilde{P})- \begin{bmatrix}0 \\K-\Tilde{K}\end{bmatrix}\right)\right)\\
\succeq & \left(\boldsymbol{\delta}^*\Tilde{P}-\begin{bmatrix}
    I\\
    \Tilde{K}
\end{bmatrix}\right)^\top\left(\boldsymbol{\delta}^*\Tilde{P}-\begin{bmatrix}
    I\\
    \Tilde{K}
\end{bmatrix}\right) \\
+& \mathcal{H}\left(\left(\boldsymbol{\delta}^*\Tilde{P}-\begin{bmatrix}
    I\\
    \Tilde{K}
\end{bmatrix}\right)^\top\left(\boldsymbol{\delta}^*(P-\Tilde{P})- \begin{bmatrix}0 \\K-\Tilde{K}\end{bmatrix}\right)\right)\\
\triangleq & L(K,P|\Tilde{K},\Tilde{P}).
\end{aligned}
\end{equation}
Thus, we have
\begin{equation}\label{bilinear term scaling}
    \mathcal{H}\left({P\boldsymbol{\delta}^*}^\top\begin{bmatrix}
    I\\
    K
\end{bmatrix}\right) \preceq \frac{1}{2}\mathcal{G}\left(\boldsymbol{\delta}^*P+\begin{bmatrix}
    I\\
    K
\end{bmatrix}\right) - \frac{1}{2}L(K,P|\Tilde{K},\Tilde{P}),
\end{equation}
which motivates a sufficient condition for \eqref{stabilizing temp2} as follows:
\begin{equation*}
    \begin{bmatrix}
        \frac{1}{2}\mathcal{G}\left(\boldsymbol{\delta}^*P+\begin{bmatrix}
    I\\
    K
\end{bmatrix}\right) - \frac{1}{2}L(K,P|\Tilde{K},\Tilde{P})  & \begin{bmatrix}
        I\\
        K
    \end{bmatrix}^\top & P\\
    \\[-10pt]
    \begin{bmatrix}
        I\\
        K
    \end{bmatrix} & -\mathbf{A}^* & 0\\
    P &0 &-I
    \end{bmatrix} \prec 0.
\end{equation*}
By applying the Schur complement to the top-left block and incorporating the structural consideration into the objective, we have the following SDP problem:
\begin{subequations}\label{stabilizing relaxed problem}
    \begin{align}
    &\min_{K,P} \quad  ||K\circ I_{S^c}||_F^2 \qquad \mathrm{s.t.}\label{stabilizing relaxed problem eq1}\\
    &\scalebox{0.94}{$\begin{bmatrix}
        - \frac{1}{2}L(K,P|\Tilde{K},\Tilde{P}) &\frac{1}{\sqrt{2}}\left(\boldsymbol{\delta}^*P+\begin{bmatrix}
    I\\
    K
\end{bmatrix}\right)^\top & \begin{bmatrix}
        I\\
        K
    \end{bmatrix}^\top &P\\
    \frac{1}{\sqrt{2}}\left(\boldsymbol{\delta}^*P+\begin{bmatrix}
    I\\
    K
\end{bmatrix}\right) & -I & 0 &0\\
    \begin{bmatrix}
        I\\
        K
    \end{bmatrix} &0 & -\mathbf{A}^* & 0\\
    P &0 &0 &-I
    \end{bmatrix} \prec 0$}\label{stabilizing relaxed problem eq2}\\
    &P \succ 0
    \end{align}
\end{subequations}
The objective function $||K\circ I_{S^c}||_F^2$ penalizes the nonzero entries of $K$ that are not allowed by the structural constraint. Problem \eqref{stabilizing relaxed problem} has several important features, which are summarized in the following theorem.
\begin{mythr}\label{stabilizing Theorem}
    Suppose that $(\Tilde{K},\Tilde{P})$ is feasible for \eqref{stabilizing temp1}. Then, Problem \eqref{stabilizing relaxed problem} is feasible and any optimal solution $(K^*,P^*)$ is feasible for \eqref{stabilizing temp1}. Moreover, $||K^*\circ I_{S^c}||_F^2 \leq ||\Tilde{K}\circ I_{S^c}||_F^2.$
\end{mythr}
\quad \textit{Proof}: We first prove the first statement. If $(\Tilde{K},\Tilde{P})$ is feasible for \eqref{stabilizing temp1}, then at least it is feasible for Problem \eqref{stabilizing relaxed problem}, since \eqref{stabilizing temp1} and \eqref{stabilizing relaxed problem eq2} are equivalent when $K=\Tilde{K}, P=\Tilde{P}$. Therefore, Problem \eqref{stabilizing relaxed problem} is feasible. Moreover, since \eqref{stabilizing relaxed problem eq2} is a sufficient condition for \eqref{stabilizing temp1}, any feasible solution to Problem \eqref{stabilizing relaxed problem} is feasible for \eqref{stabilizing temp1}. Thus, the optimal solution $(K^*,P^*)$ must be feasible for \eqref{stabilizing temp1}. For the second statement, since at least $(\Tilde{K},\Tilde{P})$ is feasible for Problem \eqref{stabilizing relaxed problem}, the optimal solution $K^*,P^*$ must satisfy $||K^*\circ I_{S^c}||_F^2 \leq ||\Tilde{K}\circ I_{S^c}||_F^2$. The proof is completed. $\hfill \square$\\

\begin{algorithm}[t]
\caption{Structured stabilizing controller design.}
\textbf{Output:} $K^*$;\\
    Solve Problem \eqref{data driven unstructured stabilizing} and denote the solution as $(\Tilde{X},\Tilde{Y})$. Initialize $P_0 = \Tilde{X}^{-1}, K_0 = \Tilde{Y}\Tilde{X}^{-1}, k=0$;\\
    \Repeat{$||K\circ I_{S^c}||_F<\epsilon_T$}{
    - Solve Problem \eqref{stabilizing relaxed problem} with $\Tilde{P} = P_k, \Tilde{K} = K_k$;\\
    - Assign the solutions to $K_{k+1}, P_{k+1}$;\\
    - $k = k + 1$;
    }
\Return $K_k$.
\label{structured stabilizing controller design algorithm}
\end{algorithm}

\cref{stabilizing Theorem} motivates an iterative algorithm to compute a structured stabilizing controller. As long as we initialize a feasible solution $(\Tilde{K},\Tilde{P})$ for \eqref{stabilizing temp1}, and use the optimal solution $(K^*,P^*)$ of \eqref{stabilizing relaxed problem} as the new $(\Tilde{K},\Tilde{P})$ to solve \eqref{stabilizing relaxed problem} again, we can ensure that this process is recursively feasible and the objective value is non-increasing. The algorithm is summarized in \cref{structured stabilizing controller design algorithm}, where $\epsilon_T$ is a small positive threshold for convergence checks.

\subsection{$H_2$ Control}
Similar to the stabilizing control case, we first replace $\begin{bmatrix}
    A &B\end{bmatrix}$ in \eqref{H2 eq1} using \cref{critical alternative}. After applying the Schur complement, we require
\begin{equation*}
    \mathcal{H}\left(P{\boldsymbol{\delta}^*}^\top \begin{bmatrix}       I\\
        K
    \end{bmatrix}\right)
    +
    \mathcal{H}\left(P\Gamma^\top ({\mathbf{A}^*}^{-\frac{1}{2}})^\top \begin{bmatrix}
        I\\
        K
    \end{bmatrix}\right) + \mathcal{G}(C_K)\prec 0
\end{equation*}
to hold for all $\|\Gamma\|\leq 1$. By the strict Petersen's lemma, this is equivalent to verifying
\begin{equation}\label{H2 temp1}
    \mathcal{H}\left(P{\boldsymbol{\delta}^*}^\top \begin{bmatrix}       I\\
        K
    \end{bmatrix}\right) + \lambda PP + \frac{1}{\lambda}\begin{bmatrix}
        I\\
        K
    \end{bmatrix}^\top {\mathbf{A}^*}^{-1} \begin{bmatrix}
        I\\
        K
    \end{bmatrix} + \mathcal{G}(C_K)\prec 0
\end{equation}
for some $\lambda>0$. Similar to \cref{Stabilizing Control}, by letting $X=P^{-1}$ and $Y=KX$, multiplying \eqref{H2 temp1} on both sides by $X$, and applying the Schur complement, we obtain the optimal unstructured $H_2$ controller via the following SDP:
\begin{subequations}\label{data driven unstructured H2}
    \begin{align}
        &\min_{\gamma,X,Y,Z,\lambda} \quad  \gamma \qquad \mathrm{s.t.}\\
        & 
    \begin{bmatrix}
        \mathcal{H}\left({\boldsymbol{\delta}^*}^\top\begin{bmatrix}
        X\\
        Y
    \end{bmatrix}\right)  + \lambda I & \begin{bmatrix}
        X\\
        Y
    \end{bmatrix}^\top & (CX+DY)^\top\\
    \\[-10pt]
        \begin{bmatrix}
        X\\
        Y
    \end{bmatrix} & -\lambda\mathbf{A}^* &0\\
    CX+DY &0 &-I
    \end{bmatrix} \prec 0\\
    &
    \begin{bmatrix}
        Z & G^\top\\
        G & X
    \end{bmatrix} \succeq 0, \mathrm{Tr}(Z) \leq \gamma^2, \lambda >0.
    \end{align}
\end{subequations}
Denote the optimal solution to \eqref{data driven unstructured H2} as $(\gamma^*,X^*,Y^*,Z^*,\lambda^*)$; the optimal $H_2$ controller can then be computed as $K^* = Y^*{X^*}^{-1}$ with the $H_2$ performance bound $\gamma^*$. Similar to the stabilizing control case, we cannot incorporate the structural constraint $K\in S$ due to the change of variable $Y=KX$. Instead, we work directly with \eqref{H2 temp1}. Applying the Schur complement, we obtain
\begin{equation}\label{H2 temp2}
    \begin{bmatrix}
        \mathcal{H}\left(P{\boldsymbol{\delta}^*}^\top\begin{bmatrix}
        I\\
        K
    \end{bmatrix}\right) & \begin{bmatrix}
        I\\
        K
    \end{bmatrix}^\top & P  & (C+DK)^\top\\
    \\[-10pt]
    \begin{bmatrix}
        I\\
        K
    \end{bmatrix} & -\lambda\mathbf{A}^* & 0 &0\\
    P &0 &-\frac{1}{\lambda}I &0\\
    C+DK &0 &0 &-I
    \end{bmatrix} \prec 0
\end{equation}
There are two challenges in handling \eqref{H2 temp2}. The first is the bilinear term in the top-left block, which can be addressed using the same technique as in the stabilizing control case. The second challenge is the non-convex (concave) term $-\frac{1}{\lambda}I$ in the third diagonal block. To tackle this, we utilize the following inequality for a given $\Tilde{\lambda}>0$:
\begin{equation}\label{inverse term scaling}
    -\frac{1}{\lambda}I \preceq -\frac{1}{\Tilde{\lambda}}I + \frac{1}{\Tilde{\lambda}^2}(\lambda - \Tilde{\lambda})I.
\end{equation}
This inequality holds due to the concavity of $-\frac{1}{\lambda}$. Applying \eqref{bilinear term scaling}, \eqref{inverse term scaling}, and the Schur complement, we can now convexify \eqref{H2 temp2} to \eqref{H2 relaxed problem eq2}. It is easy to verify that \eqref{H2 relaxed problem eq2} is a sufficient condition for \eqref{H2 temp2}. In the objective function \eqref{H2 relaxed problem eq1}, we minimize the $H_2$ performance bound $\gamma$ plus a regularization term $\beta||K\circ I_{S^c}||_F^2$ to promote the structural constraint $K\in S$. The parameter $\beta>0$ is a weight reflecting the emphasis on the structural constraint. The full optimization problem is shown in Problem~\eqref{H2 relaxed problem}.
\begin{subequations}\label{H2 relaxed problem}
    \begin{align}
    &\min_{\gamma,K,P,Z,\lambda} \quad  \gamma + \beta||K\circ I_{S^c}||_F^2 \qquad \mathrm{s.t.} \label{H2 relaxed problem eq1}\\
    &\begin{aligned}&\begin{bmatrix}
        - \frac{1}{2}L(K,P|\Tilde{K},\Tilde{P}) &* & *&*  &*\\
    \frac{1}{\sqrt{2}}\left(\boldsymbol{\delta}^*P+\begin{bmatrix}
    I\\
    K
\end{bmatrix}\right) & -I &* &* &*\\
\\[-10pt]
    \begin{bmatrix}
        I\\
        K
    \end{bmatrix} &0 & -\lambda \mathbf{A}^* & * &*\\
    P &0 &0 &\Lambda &*\\
    C+DK &0 &0 &0 &-I
    \end{bmatrix} \prec 0 \label{H2 relaxed problem eq2}\\
    &\Lambda = \left(-1/\Tilde{\lambda}+(1/\Tilde{\lambda}^2)(\lambda-\Tilde{\lambda})\right) I,
    \end{aligned}\\
    &P \succ 0, Z-G^\top PG \succeq 0, \mathrm{Tr}(Z) \leq \gamma^2, \lambda > 0.\label{H2 relaxed problem eq3}
    \end{align}
\end{subequations}
Some important features of Problem \eqref{H2 relaxed problem} are summarized in the following theorem.

\begin{algorithm}[t]
\caption{Structured $H_2(H_\infty)$ controller design.}
\textbf{Output:} $K^*,\gamma^*$;\\
    $\text{(I)}=\begin{cases}
    \text{Problem \eqref{data driven unstructured H2},}  & \text{Unstructured}\ H_2\ \text{controller design}\\
    \text{Problem \eqref{data-driven unstructured hinf},} & \text{Unstructured}\ H_\infty\ \text{controller design}\\
    \end{cases}$\\
    $\text{(P)}=\begin{cases}
    \text{Problem \eqref{H2 relaxed problem},}  & \text{Structured}\ H_2\ \text{controller design}\\
    \text{Problem \eqref{hinf relaxed problem},} & \text{Structured}\ H_\infty\ \text{controller design}\\
    \end{cases}$\\
    Solve Problem (I) and denote the solution as $(\Tilde{X},\Tilde{Y},\Tilde{\lambda})$. Initialize $P_0 = \Tilde{X}^{-1}, K_0 = \Tilde{Y}\Tilde{X}^{-1}, \lambda_0 = \Tilde{\lambda}, \beta = 1, \mu > 1, k=0$;\\
    \Repeat{$||P_{k}-P_{k-1}||_F<\epsilon_T$ and $||K_{k}-K_{k-1}||_F<\epsilon_T$}{
    - Solve Problem (P) with $\Tilde{P} = P_k, \Tilde{K} = K_k, \Tilde{\lambda} = \lambda_k$;\\
    - Assign the solutions to $K_{k+1}, P_{k+1}, \lambda_{k+1}, \gamma_{k+1}$;\\
    - \textbf{if} $\beta<1e6$ \textbf{then} $\beta = \mu\beta$ \textbf{end if}\\
    - $k = k + 1$;
    }
\Return $K_k,\gamma_k$.
\label{structured h2(hinf) controller design algorithm}
\end{algorithm}

\begin{mythr}\label{h2 Theorem}
    Suppose that $(\Tilde{K},\Tilde{P},\Tilde{\lambda})$ is feasible for \eqref{H2 temp1}. Then, Problem \eqref{H2 relaxed problem} is feasible and any optimal solution $(K^*,P^*,\lambda^*)$ is feasible for \eqref{H2 temp1}. 
\end{mythr}
\quad \textit{Proof}: Similar to \cref{stabilizing Theorem}, if $(\Tilde{K},\Tilde{P},\Tilde{\lambda})$ is feasible for \eqref{H2 temp1}, then at least $(\Tilde{K},\Tilde{P},\Tilde{\lambda})$ is feasible for \eqref{H2 relaxed problem} since \eqref{H2 temp1} and \eqref{H2 relaxed problem eq2} are equivalent when $K=\Tilde{K}, P=\Tilde{P}$, which shows Problem \eqref{H2 relaxed problem} is feasible. Since \eqref{H2 relaxed problem eq2} is a sufficient condition for \eqref{H2 temp1}, any feasible solution to Problem \eqref{H2 relaxed problem}, including the optimal solution $(K^*,P^*,\lambda^*)$, is feasible for \eqref{H2 temp1}. The proof is completed. $\hfill \square$\\
Similar to the stabilizing control case, we can develop an iterative algorithm to compute the optimal structured $H_2$ controller based on \cref{h2 Theorem}. This process remains recursively feasible provided that it is initialized with a feasible solution $(\Tilde{K},\Tilde{P},\Tilde{\lambda})$ for \eqref{H2 temp1}. In practice, the weight $\beta$ is gradually increased to enforce the structural constraint. The algorithm is presented in the $H_2$ case of \cref{structured h2(hinf) controller design algorithm}.

\subsection{$H_\infty$ Control}
The top-left block of \eqref{H2 eq1} and \eqref{Hinf eq} are the same. Therefore, the technique used in the $H_2$ control case can be directly extended to the $H_\infty$ control case. For brevity, we omit the derivation details and present the unstructured and structured $H_\infty$ controller synthesis problems as Problems \eqref{data-driven unstructured hinf} and \eqref{hinf relaxed problem}.

\begin{subequations}\label{data-driven unstructured hinf}
    \begin{align}
        &\min_{\gamma,X,Y,\lambda} \quad  \gamma \qquad \mathrm{s.t.}\\
        & \scalebox{0.98}{$
    \begin{bmatrix}
        \mathcal{H}\left({\boldsymbol{\delta}^*}^\top\begin{bmatrix}
        X\\
        Y
    \end{bmatrix}\right)  + \lambda I & \begin{bmatrix}
        X\\
        Y
    \end{bmatrix}^\top &G & (CX+DY)^\top\\
    \\[-10pt]
        \begin{bmatrix}
        X\\
        Y
    \end{bmatrix} & -\lambda\mathbf{A}^* &0 &0\\
    \\[-10pt]
    G^\top &0 &-\gamma I &H^\top\\
    CX+DY &0 &H &-\gamma I
    \end{bmatrix} \prec 0$}\\
    &X \succ 0, \lambda >0.
    \end{align}
\end{subequations}

\begin{subequations}\label{hinf relaxed problem}
    \begin{align}
    &\min_{\gamma,K,P,\lambda} \quad  \gamma + \beta||K\circ I_{S^c}||_F^2 \qquad \mathrm{s.t.}\\
    &\begin{aligned}&\begin{bmatrix}
        - \frac{1}{2}L(K,P|\Tilde{K},\Tilde{P}) &* & * &*  &* &*\\
    \frac{1}{\sqrt{2}}\left(\boldsymbol{\delta}^*P+\begin{bmatrix}
    I\\
    K
\end{bmatrix}\right) & -I & * &* &* &*\\
\\[-10pt]
    \begin{bmatrix}
        I\\
        K
    \end{bmatrix} &0 & -\lambda\mathbf{A}^* & * &* &*\\
    P &0 &0 &\Lambda &* &*\\
    G^\top P &0 &0 &0 &-\gamma I &*\\
    C+DK &0 &0 &0  &H &-\gamma I
    \end{bmatrix} \prec 0\\
    &\Lambda = \left(-1/\Tilde{\lambda}+(1/\Tilde{\lambda}^2)(\lambda-\Tilde{\lambda})\right) I,
    \end{aligned}\\
    &P \succ 0, \lambda > 0.
    \end{align}
\end{subequations}

Since the analysis process is similar to the $H_2$ control case, \cref{h2 Theorem} can be directly extended to the $H_\infty$ control case. The algorithm is summarized in the $H_\infty$ case of \cref{structured h2(hinf) controller design algorithm}.

\begin{myrem}
    In this paper, we do not provide convergence guarantees, as there may not exist a stabilizing structured controller for the ground truth system, and verifying its existence is NP-hard \cite{blondel1997np}. Nevertheless, under reasonable structural constraints, our algorithms often succeed in finding a structured controller with satisfactory performance, even in cases where the method in \cite{jovanovic2016controller} fails to produce feasible solutions. Therefore, our approach offers a promising direction for data-driven structured controller design problems.
\end{myrem}

\section{Simulations}\label{Simulations}
In this section, we provide simulation results to demonstrate the effectiveness of our proposed methods. Specifically, we compare our methods with the existing approach that fixes $X$ to be diagonal and imposes the structural constraint on $Y$ \cite{jovanovic2016controller}. Throughout this section, we set $T_s = 0.1$, $\mu=2$, and $\epsilon_T = 0.01$. In the result tables, the $(A,B)$ column denotes the model-based case, while the other columns represent the data-driven cases with different noise bounds or data lengths.

\begin{table}[!t]
\centering
\caption{$H_2$ performance bound comparison for $T=100$.}
\begin{tabular}{cccccc}
\toprule
& {Design}              & $(A,B)$ & $\epsilon=0.01$ & $\epsilon=0.03$  &$\epsilon=0.05$\\
\midrule
& $X$ diag      &Infeasible &Infeasible &Infeasible &Infeasible\\
& \textbf{Ours}     & 2.2831 & 2.5103 & 3.1647  & 4.4266\\
\bottomrule
\end{tabular}
\label{H2 comparison for fixed T}
\end{table}

\begin{table}[!t]
\centering
\caption{$H_2$ performance bound comparison for $\epsilon=0.01$.}
\begin{tabular}{cccccc}
\toprule
& {Design}              & $(A,B)$ & $T=60$ & $T=80$  &$T=100$\\
\midrule
& $X$ diag      &Infeasible &Infeasible &Infeasible &Infeasible\\
& \textbf{Ours}     & 2.2831 & 2.8976 & 2.6062  & 2.5129\\
\bottomrule
\end{tabular}
\label{H2 comparison for fixed noise bound}
\end{table}

\begin{table}[!t]
\centering
\caption{$H_\infty$ performance bound comparison for $T=100$.}
\begin{tabular}{cccccc}
\toprule
& {Design}              & $(A,B)$ & $\epsilon=0.01$ & $\epsilon=0.03$  &$\epsilon=0.05$\\
\midrule
& $X$ diag      &Infeasible &Infeasible &Infeasible &Infeasible\\
& \textbf{Ours}     & 1.7533 & 1.8501 & 2.1117  & 2.5398\\
\bottomrule
\end{tabular}
\label{Hinf comparison for fixed T}
\end{table}

\begin{table}[!t]
\centering
\caption{$H_\infty$ performance bound comparison for $\epsilon=0.01$.}
\begin{tabular}{cccccc}
\toprule
& {Design}              & $(A,B)$ & $T=60$ & $T=80$  &$T=100$\\
\midrule
& $X$ diag      &Infeasible &Infeasible &Infeasible &Infeasible\\
& \textbf{Ours}     & 1.7533 & 1.8900 & 1.8619  & 1.8479\\
\bottomrule
\end{tabular}
\label{Hinf comparison for fixed noise bound}
\end{table}

Consider a mass-spring system with two masses connected in a line \cite{lin2011augmented}. The system matrices are given by
\begin{equation*}
    A = \begin{bmatrix}
    0 &I\\
    T &0
    \end{bmatrix}, \quad B = \begin{bmatrix}
    0\\I
    \end{bmatrix}, \quad G = \begin{bmatrix}
    0\\I
    \end{bmatrix},
\end{equation*}
where $T\in \mathbb{R}^{2\times2}$ is a Toeplitz matrix with $-2$ on its main diagonal and $1$ on its first upper and lower diagonals. We first consider \cref{problem 1}. We impose a structural constraint with
\begin{equation*}
    I_S = \begin{bmatrix}
    0 &1 &1 &0\\
    0 &1 &1 &0
    \end{bmatrix}.
\end{equation*}
We collect data with $T=100$ and $\epsilon=0.01$. The method proposed in \cite{jovanovic2016controller}, which requires $X$ to be diagonal and $Y\circ I_{S^c}=0$ in Problem \eqref{data driven unstructured stabilizing}, leads to infeasibility. In contrast, after running \cref{structured stabilizing controller design algorithm}, we obtain
\begin{equation*}
    K = \begin{bmatrix}
    0 & 0.9078 & -2.0189 & 0\\
    0 & 0.3254 & 0.3022 & 0
    \end{bmatrix},
\end{equation*}
which renders the ground truth closed-loop system stable. Next, we consider \cref{problem 2}. We choose 
\begin{equation*}
    C = \begin{bmatrix}
        I \\ 0
        \end{bmatrix}, \quad
    D = \begin{bmatrix}
        0 \\ I
        \end{bmatrix}, \quad
    I_S = \begin{bmatrix}
    1 &1 &0 &0\\
    0 &0 &1 &1
    \end{bmatrix}.
\end{equation*}
The simulation results are reported in \cref{H2 comparison for fixed T} and \cref{H2 comparison for fixed noise bound}. The results show that the method in \cite{jovanovic2016controller} leads to infeasibility in all cases, while our method remains feasible. In addition, it is worth noting that as the noise bound $\epsilon$ increases or the data length $T$ decreases, the $H_2$ performance bound becomes worse, which is consistent with intuition. We finally consider \cref{problem 3}. We choose 
\begin{equation*}
        C = \begin{bmatrix}
        I \\ 0
        \end{bmatrix}, \quad
    D = \begin{bmatrix}
        0 \\ I
        \end{bmatrix}, \quad H = \begin{bmatrix}
        0 \\ I
        \end{bmatrix}, \quad
    I_S = \begin{bmatrix}
    0 &0 &1 &1\\
    1 &0 &0 &0
    \end{bmatrix}.
\end{equation*}
The simulation results are reported in \cref{Hinf comparison for fixed T} and \cref{Hinf comparison for fixed noise bound}. The results are similar to the $H_2$ case, further demonstrating the effectiveness of our proposed method.

\section{Conclusion}\label{Conclusion}
In this paper, we address the problem of data-driven structured controller design for continuous-time LTI systems. Three control objectives are considered, including stabilization, $H_2$ control, and $H_\infty$ control. We first describe the minimal matrix ellipsoid containing all possible system matrices that can generate the collected data. Then, for each control objective, we derive the data-driven unstructured controller design method via an SDP. To handle the structural constraint, we propose novel linearization techniques to relax the original nonconvex problems into SDPs. Based on these relaxations, we develop iterative algorithms to compute structured controllers. Simulation results demonstrate the effectiveness of our methods.

\footnotesize{
\bibliographystyle{IEEEtran}
\bibliography{reference}
}
\end{document}